\documentclass{amsart}

\usepackage{amscd}
\usepackage{amssymb}
\theoremstyle{plain}

\newtheorem{thm}{Theorem}

\newtheorem{lem}[thm]{Lemma}

\theoremstyle{definition}
\newtheorem{df}[thm]{Definition}
\newtheorem{rem}[thm]{Remark}

\newtheorem{calc}[thm]{Calculation}

\author{Michael Robinson} 

\title{IMEX method convergence for a parabolic equation}

\begin{document}
\begin{abstract}
Although implicit-explicit (IMEX) methods for approximating solutions to
semilinear parabolic equations are relatively standard, most recent
works examine the case of a fully discretized model.  We show that by
discretizing time only, one can obtain an elementary convergence
result for an implicit-explicit method.  This convergence result is
strong enough to imply existence and uniqueness of solutions to a
class of semilinear parabolic equations.
\end{abstract}


\maketitle

\section{Introduction}

The use of implicit-explicit (IMEX) methods for approximating
semilinear parabolic equations is well-established
\cite{AscherRuuthWetton}.  Many of the recent works on these methods
employ discretizations in both space and time.  These fully discrete
approximations can be computed directly by a computer.  However, one
can obtain a stronger condition for convergence of the approximation
if only the time dimension is discretized \cite{Crouzeix}.  We show
how an even stronger condition for convergence is met by the Cauchy
problem for
\begin{equation}
\label{pde2}
\frac{\partial u(x,t)}{\partial t}=\Delta u(x,t) + \sum_{i=0}^\infty
a_i(x) u^i(x,t),
\end{equation}
where $a_i \in L^1(\mathbb{R}^n) \cap L^\infty(\mathbb{R}^n)$, and
how convergence of this method provides an elementary proof of
existence and uniqueness of solutions.  Existence and uniqueness of
solutions for \eqref{pde2} under reasonable initial conditions have
been known for some time.  For instance, \cite{Henry} and
\cite{ZeidlerIIA} contain straightforward proofs using semigroup
methods.  The purpose of this paper is to show how a {\it more
elementary} proof can be obtained from a sequence of explicitly
computed discrete-time approximations.

The Cauchy problem for \eqref{pde2} arises in a variety of settings.
Notably, some reaction-diffusion equations are of this
form \cite{FiedlerScheel}.  Another application is the special case
\begin{equation*}
\frac{\partial u(x,t)}{\partial t}=\Delta u(x,t) - u^2(x,t) + a_0(x),
\end{equation*}
where $a_0$ is a nonzero function of $x$.  This situation corresponds
to a spatially-dependent logistic equation with a diffusion term,
which can be thought of as a toy model of population growth with
migration.

Following \cite{Crouzeix}, the approximation to be used is
\begin{equation}
\label{e_i_deriv}
u_{n+1}=(I-h \Delta)^{-1}(u_n + h \sum_{i=0}^\infty a_i u_n^i),
\end{equation}
which is obtained by inverting the linear portion of a discrete
version of \eqref{pde2}.  For brevity, we shall call \eqref{e_i_deriv}
{\it the} implicit-explicit method.  (In the summary paper
\cite{AscherRuuthWetton}, this is called an SBDF method, to
distinguish it from other implicit-explicit methods.)  One can compute
the operator $(I-h\Delta)^{-1}$ explicitly using Fourier transform
methods, and obtain a proof of the numerical stability of the
iteration as a whole.

\section{A version of the fundamental inequality}

In order to simplify the algebraic expressions, we make the following definitions.

\begin{df}
Let 
\begin{equation}
\label{def_F}
F(u(x,t))=\Delta u(x,t) + \sum_{i=0}^\infty a_i(x) u^i(x,t),
\end{equation}
and
\begin{equation}
\label{def_G}
G(u(x,t))=\sum_{i=0}^\infty a_i(x) u^i(x,t).
\end{equation}
\end{df}
\begin{df}
\label{def_g}
Define the analytic functions
\begin{equation}
\label{def_g_1}
g_1(z)=\sum_{i=0}^\infty \|a_i\|_1 z^i,
\end{equation}
and
\begin{equation}
\label{def_g_infty}
g_\infty(z)=\sum_{i=0}^\infty \|a_i\|_\infty z^i.
\end{equation}
\end{df}

Since we do not discretize the spatial dimension, we can employ some
of the theory of ordinary differential equations.  We therefore first
prove a variant of the fundamental inequality for \eqref{pde2} as is
done in \cite{HubbardWest}.  The fundamental inequality gives a
sufficient condition for approximate solutions to converge.  A
slightly weaker version of Lemma \ref{excond_lem} was obtained in
Theorem 3.1 of \cite{Crouzeix}, where the existence of solutions was
required.

\begin{lem}
\label{excond_lem}
Suppose $\{u_i\}_{i=1}^\infty$ is a sequence of piecewise $C^1$
functions $u_i:[0,T] \rightarrow C^2(\mathbb{R}^n)\cap
L^1(\mathbb{R}^n) \cap L^\infty(\mathbb{R}^n)$, such that 
\begin{enumerate}
\item there exist $A,B>0$ so that for each $i$ and $t \in [0,T]$,
  $\|u_i(t)\|_1 \le A$ and $\|u_i(t)\|_\infty \le B$,
\item for each $i$ and $t\in [0,T]$, the series
  $g_1(\|u_i(t)\|_1)$ and $g_\infty(\|u_i(t)\|_\infty)$ converge,
\item for each $t \in [0,T]$,
  $\|\frac{d}{dt} u_i(t) - F(u_i(t)) \|_\infty < \epsilon_i$ and
  $\lim_{i\rightarrow\infty}\epsilon_i = 0$, and
\item $u_1(0) = u_i(0)$ for all $i \ge 0$
\end{enumerate}
Then for each $t \in [0,T]$, $\{u_i(t)\}_{i=1}^\infty$ is a Cauchy
sequence in $L^2(\mathbb{R}^n)$.
\begin{proof}
Let $i,j>0$ be given.  Let $\eta(t)=\|u_i(t) -
u_j(t)\|_2^2=\int{(u_i(t)-u_j(t))^2 dx}$.  Notice that the fourth
condition in the hypothesis gives $\eta(0)=0$.  
\begin{equation*}
\eta'(t)=2\int{\left(u_i'(t)-u_j'(t)\right)\left(u_i(t)-u_j(t)\right) dx}.
\end{equation*}
But, $\|\frac{d}{dt} u_i(t) - F(u_i(t)) \|_\infty < \epsilon_i$ is
equivalent to the statement that for each $t \in [0,T]$ and $x \in
\mathbb{R}^n$,  
\begin{equation*}
F(u_i(x,t)) - \epsilon_i < u_i'(x,t) < F(u_i(x,t)) + \epsilon_i,
\end{equation*}
giving
\begin{eqnarray*}
\eta'(t) &\le& 2 \int{\left(F(u_i(t))-F(u_j(t)) \right)(u_i(t)-u_j(t))
  dx} + 2(\epsilon_i + \epsilon_j)\int{|u_i(t)-u_j(t)| dx}\\
&\le&2 \int{\left(\Delta u_i(t) + G(u_i(t)) - \Delta u_j(t) -
  G(u_j(t))\right)(u_i(t) - u_j(t)) dx} \\
 &&+ 2(\epsilon_i + \epsilon_j)\|u_i(t)-u_j(t)\|_1 \\
&\le&2 \int{\left(\Delta(u_i(t)-u_j(t))\right)(u_i(t)-u_j(t)) dx} \\
 &&+ 2 \int{\left(G(u_i(t)) - G(u_j(t))\right)(u_i(t)-u_j(t)) dx} 
 + 2(\epsilon_i + \epsilon_j)\|u_i(t)-u_j(t)\|_1 \\
&\le& -2 \int{\|\nabla(u_i(t)-u_j(t))\|^2 dx} + 2
  \|G(u_i(t)) - G(u_j(t))\|_2 \|u_i(t)-u_j(t)\|_2 \\
 &&+ 2(\epsilon_i + \epsilon_j)\|u_i(t)-u_j(t)\|_1 \\
&\le& 2 \|G(u_i(t)) - G(u_j(t))\|_2 \|u_i(t)-u_j(t)\|_2 +
 2(\epsilon_i + \epsilon_j)\|u_i(t)-u_j(t)\|_1. \\
\end{eqnarray*}
Now also 
\begin{eqnarray*}
\|G(u_i(t)) &-& G(u_j(t))\|_2 = \left \|\sum_{k=0}^\infty a_k
(u_i^k(t)-u_j^k(t)) \right \|_2\\
&\le&\sum_{k=0}^\infty \|a_k\|_\infty \left \|u_i^k(t)-u_j^k(t) \right
\|_2 \\
&\le&\sum_{k=0}^\infty \|a_k\|_\infty
\sqrt{\int{\left(u_i^k(x,t)-u_j^k(x,t)\right)^2 dx}} \\
&\le&\sum_{k=0}^\infty \|a_k\|_\infty
\sqrt{\int{\left(u_i(x,t)-u_j(x,t)\right)^2 \left(\sum_{m=0}^{k-1} u_i^m(x,t)
    u_j^{k-m-1}(x,t)\right)^2 dx}} \\
&\le&\sum_{k=0}^\infty \|a_k\|_\infty \left \| \sum_{m=0}^{k-1} u_i^m(t)
    u_j^{k-m-1}(t) \right \|_\infty \|u_i(t)-u_j(t)\|_2 \\
&\le&\left( \sum_{k=0}^\infty \|a_k\|_\infty k B^{k-1}
    \right)\|u_i(t)-u_j(t)\|_2 \\
&\le&g_\infty'(B) \|u_i(t)-u_j(t)\|_2, \\
\end{eqnarray*}
which allows
\begin{eqnarray*}
\eta'(t) &\le& 2 g_\infty'(B) \|u_i(t)-u_j(t)\|_2^2 +  2(\epsilon_i +
\epsilon_j)\|u_i(t)-u_j(t)\|_1. \\
&\le&2 g_\infty'(B) \eta(t) +  2(\epsilon_i +
\epsilon_j)\|u_i(t)-u_j(t)\|_1. \\
\end{eqnarray*}
\begin{eqnarray*}
\eta'(t)-2 g_\infty'(B) \eta(t) &\le& 2(\epsilon_i +
\epsilon_j)\|u_i(t)-u_j(t)\|_1 \\
\frac{d}{dt}\left(\eta(t) e^{-2 g_\infty'(B) t } \right) 
&\le& 2(\epsilon_i +\epsilon_j) e^{-2 g_\infty'(B) t}\|u_i(t)-u_j(t)\|_1, \\
\end{eqnarray*}
so (recall $\eta(0)=0$)
\begin{eqnarray*}
\eta(t) &\le& \left[ 2(\epsilon_i +\epsilon_j) \int_0^t e^{-2
    g_\infty'(B) s } \|u_i(s)-u_j(s)\|_1 ds
    \right]e^{2 g_\infty'(B) t} \\
&\le& \left[ 2(\epsilon_i +\epsilon_j) \int_0^t{\|u_i(s)-u_j(s)\|_1
    ds} \right] e^{2 g_\infty'(B) t}.\\ 
&\le& 4(\epsilon_i +\epsilon_j) A t e^{2 g_\infty'(B) t}.\\
\end{eqnarray*}
Hence as $i,j \rightarrow \infty$, $\eta(t) \rightarrow 0$ for each
$t$.  Thus for each $t$, $\{u_i(t)\}_{i=1}^\infty$ is a Cauchy
sequence in $L^2(\mathbb{R}^n)$.
\end{proof}
\end{lem}

\begin{rem}
Since $C^2(\mathbb{R}^n) \cap L^1(\mathbb{R}^n) \cap L^\infty(\mathbb{R}^n)
\subseteq L^2(\mathbb{R}^n) $ and $L^2$ is complete, Lemma \ref{excond_lem}
gives conditions for existence and uniqueness of a short-time solution
to \eqref{pde2}.
\end{rem}

\begin{lem}
\label{derivative_lem}
Suppose $\{u_i(t)\}_{i=1}^\infty$ is the sequence of functions defined
in Lemma \ref{excond_lem}, and that $u(t)=\lim_{i\rightarrow\infty}
u_i(t)$ in $L^2(\mathbb{R}^n)$.  Then
\begin{equation}
\label{derivativelimit_eqn}
u'(t,x)=\lim_{i\rightarrow\infty}u_i'(t,x) \text{ for almost every $x$},
\end{equation}
wherever the limit exists.
\begin{proof}
Notice that since each $u_i(t) \in L^\infty(\mathbb{R}^n)$ and
$\|u_i(t)\|_\infty \le B$, the dominated convergence theorem allows for
each $x \in \mathbb{R}^n$
\begin{eqnarray*}
\int_0^t{\lim_{i\rightarrow\infty} u_i'(\tau,x) d\tau} &=&
  \lim_{i\rightarrow\infty} \int_0^t{u_i'(\tau,x) d\tau}\\
&=& \lim_{i\rightarrow\infty}(u_i(t,x) - u_i(0,x))\\
&=& u(t,x)-u(0,x) \text{ for almost every $x$}. 
\end{eqnarray*} 
Hence, by differentiating in $t$, 
\begin{equation*}
u'(t,x)=\lim_{i\rightarrow\infty}u_i'(t,x) \text{ for almost every $x$}.
\end{equation*}
\end{proof}
\end{lem}

\section{The implicit-explicit approximation}

In this section, we consider the case of a 1-dimensional
spatial domain, that is, $x \in \mathbb{R}$.  There is no obstruction
to extending any of these results to higher dimensions, though it
complicates the exposition unnecessarily.

As is usual, the first task is to define the function spaces to be
used.  Initial conditions will be drawn from a subspace of
$L^1(\mathbb{R})\cap L^\infty(\mathbb{R})$, as suggested by Lemma
\ref{excond_lem}, and the first four spatial derivatives will be
prescribed, for use in Lemma \ref{hubbardbnd_lem}.

\begin{df}
\label{spaces_df}
Let
\begin{equation*}
W=L^1(\mathbb{R}) \cap
L^\infty(\mathbb{R}) \cap \{f \in C^\infty(\mathbb{R}) | f
\text{ has bounded partial derivatives up to fourth order}\}.
\end{equation*}
For the remainder of this paper, we consider the case where each of
the coefficients $a_i \in W$.  Then let $X=\{f\in W | g_1(\|f\|_1) <
\infty\text{ and } g_\infty(\|f\|_\infty) < \infty\}$.  We consider
the case where the initial condition is drawn from $X$.
\end{df}

An approximate solution given by the implicit-explicit iteration will
be the piecewise linear interpolation through the iterates computed by
\eqref{e_i_deriv}.  A smoother approximation will prove to be
unnecessary, as will be shown in Lemma \ref{existunique_h_thm}.

\begin{df}
\label{hubbard_df}
Suppose $f_0$ and $h > 0$ are given.  Put 
\begin{equation}
\label{hubbard_df_pts}
f_{n+1}=(I-h \Delta)^{-1} ( f_n + h G(f_n) ).
\end{equation}
The function 
\begin{equation}
\label{hubbard_df_fcn}
u(t)=\left(1-\left(\frac{t}{h}-n(t)\right)\right) f_{n(t)} +
\left(\frac{t}{h} - n(t)\right ) f_{n(t)+1},
\end{equation}
where $n(t)=\lfloor \frac{t}{h} \rfloor$, is called the {\bf
  implicit-explicit iteration of size $h$ beginning at $f_0$}.
\end{df}

\begin{calc}
\label{invert_calc}
We explicitly compute the operator $(I-h\Delta)^{-1}$ using Fourier
transforms.  Suppose
\begin{equation*}
(I-h\Delta)u(x)=u(x)-h\Delta u(x)=f(x).
\end{equation*}
Taking the Fourier transform (with transformed variable $\omega$)
gives
\begin{equation*}
\hat{u}(\omega)+h \omega^2 \hat{u}(\omega)=\hat{f}(\omega),
\end{equation*}
\begin{equation*}
\hat{u}(\omega)=\frac{\hat{f}(\omega)}{1+h\omega^2}.
\end{equation*}
The Fourier inversion theorem yields
\begin{eqnarray*}
u(x)&=&\frac{1}{2 \pi} \int{\frac{e^{i\omega x}}{1+h\omega^2} \int{f(y)
    e^{-i\omega y} dy} d\omega}\\
&=&\int{f(y)\left(\frac{1}{2\pi}\int{\frac{e^{i\omega(x-y)}}{1+h\omega^2}
    d\omega} \right) dy}.\\
\end{eqnarray*}

Using the method of residues, this can be simplified to give
\begin{equation}
\label{explicit_op_eqn}
u(x)=\left((I-h\Delta)^{-1}f\right)(x) = \frac{1}{2\sqrt{h}}\int{ f(y)
  e^{-|y-x|/\sqrt{h}} dy }.
\end{equation}
\end{calc}

\begin{calc}
\label{bound_calc}
Bounds on the $L^1$ and $L^\infty$ operator norms of
$(I-h\Delta)^{-1}$ are now computed.  First, let $f\in
L^\infty(\mathbb{R})$.  Then
\begin{eqnarray*}
|\left((I-h\Delta)^{-1}f\right)(x)| &=& \left|\frac{1}{2\sqrt{h}}\int{ f(y)
  e^{-|y-x|/\sqrt{h}} dy }\right|\\
&\le&\|f\|_\infty \frac{1}{2\sqrt{h}} \int{e^{-|y-x|/\sqrt{h}} dy }\\
&\le&\|f\|_\infty \frac{1}{\sqrt{h}} \int_0^\infty{e^{-s/\sqrt{h}}
  ds }\\
&\le&\|f\|_\infty,
\end{eqnarray*}
so $\|(I-h\Delta)^{-1}\|_\infty \le 1$.

Now, let $f\in L^1(\mathbb{R})$.  So then
\begin{eqnarray*}
\|(I-h\Delta)^{-1}f\|_1 &=& \int_{-\infty}^\infty{ \left|
  \frac{1}{2\sqrt{h}} \int_{-\infty}^\infty{f(y)e^{-|y-x|/\sqrt{h}} dy}
  \right| dx}\\ 
&\le& \frac{1}{2\sqrt{h}} \int_{-\infty}^\infty{
  \int_{-\infty}^\infty{ |f(y)|e^{-|y-x|/\sqrt{h}} dy} dx}\\
&\le& \frac{1}{\sqrt{h}} \int_{-\infty}^\infty{|f(y)|
  \int_{0}^\infty{ e^{-|y-x|/\sqrt{h}} dx} dy}\\
&\le& \int_{-\infty}^\infty{ |f(y)| dy} = \|f\|_1,
\end{eqnarray*}
which means $\|(I-h\Delta)^{-1}\|_1 \le 1$. 
\end{calc}

The third condition of Lemma \ref{excond_lem} is a control on the
slope error of the approximation.  A bound on this error may be
established for the implicit-explicit iteration as follows.

\begin{lem}
\label{hubbardbnd_lem}
Suppose $f_0 \in X$, $h>0$.  Put $f(x,t)=
f_0(x) + t D(x)$, where
\begin{equation*}
D=\frac{(I-h\Delta)^{-1}(f_0+hG(f_0)) - f_0}{h}
\end{equation*}
Then for every $0<t<h$, 
\begin{equation}
\label{hubbardbnd_eqn}
\|f'(t)-F(f(t))\|_\infty = O(h).
\end{equation}
\begin{proof}
Recall every function in $X$ will have bounded
partial derivatives up to fourth order from Definition \ref{spaces_df}.
\begin{eqnarray*}
\|f'(t)-F(f(t))\|_\infty &=& \|D-(\Delta (f_0+tD)+G(f_0+tD))\|_\infty\\
&=&\left\|D-\left(\Delta (f_0+tD) + \sum_{i=0}^\infty a_i (f_0+tD)^i\right)
\right\|_\infty\\
&\le&\left \|D-\Delta f_0 -t \Delta D - \sum_{i=0}^\infty a_i
 \left( \sum_{j=0}^i \binom{i}{j} f_0^j(tD)^{i-j} \right) \right\|_\infty\\
&\le&\left \|D-\Delta f_0 - t \Delta D - \sum_{i=0}^\infty a_i f_0^i
 \right\|_\infty + O(h)\\
&\le&\left \|\frac{(I-h\Delta)^{-1}-I}{h}f_0-\Delta f_0\right.\\
&&\left.+((I-h\Delta)^{-1}-I)G(f_0)\right \|_\infty+O(h)\\
\end{eqnarray*}
Now, using the fact that $(I-h\Delta)^{-1}-I=(I-h\Delta)^{-1}(h\Delta)$,
\begin{eqnarray*}
\|f'(t)-F(f(t))\|_\infty
&\le&\left \|(I-h\Delta)^{-1}\Delta f_0-\Delta
f_0\right.\\
&&\left.+(I-h\Delta)^{-1}(h\Delta)G(f_0)\right \|_\infty+O(h)\\
&\le&\|(I-h\Delta)^{-1}(h\Delta)(\Delta f_0 + G(f_0))\|_\infty+O(h)\\
&\le&h \|(I-h\Delta)^{-1}(\Delta F(f_0))\|_\infty+O(h)\\
&\le&h \|(I-h\Delta)^{-1}\|_\infty\|(\Delta F(f_0))\|_\infty+O(h)=O(h)\\
\end{eqnarray*}
\end{proof}
\end{lem}

\begin{lem}
\label{existunique_h_thm}
Suppose $0<h_i \rightarrow 0$.  Let $u_i$ be the implicit-explicit
iteration of size $h_i$ beginning at $f_0 \in X$ on $t \in [0,T]$.
Then provided there exist $A,B>0$ such that for each $i$ and $t \in
[0,T]$, $\|u_i(t)\|_1 \le A$ and $\|u_i(t)\|_\infty \le B$, then the
sequence $\{u_i(t)\}_{i=1}^\infty$ converges pointwise to a function
in $t$.  The limit function is piecewise differentiable in $t$.
\begin{proof}
Let $u_i$ be the implicit-explicit iteration of size $h_i$.  By Lemma
\ref{hubbardbnd_lem}, the slope error is bounded:
\begin{equation*}
\|u_i'(t)-F(u_i(t))\|_\infty = O(h_i) = \epsilon_i.
\end{equation*}
Notice that $\epsilon_i \rightarrow 0$.  Then, since $X \subset
C^2(\mathbb{R}^n)$, Lemma \ref{excond_lem} applies, giving a pointwise
limit function $u(t)$.  Finally, since the slope error uniformly
vanishes, Lemma \ref{derivative_lem} implies that the solution is
piecewise differentiable.
\end{proof}
\end{lem}

\section{``{\it A priori} estimates'' for the approximate solutions}

Now we demonstrate that the implicit-explicit method converges for all
initial conditions in $X$.  Specifically, for each $f_0 \in X$, there
exist $A,B>0$ such that for each $i$ and $t \in [0,T]$, $\|u_i(t)\|_1
\le A$ and $\|u_i(t)\|_\infty \le B$, given sufficiently small $T$.
We begin by recalling that from Calculation \ref{bound_calc}, the
$L^\infty$-norm of $(I-h\Delta)^{-1}$ is less than one.  This means
that for the implicit-explicit iteration,
\begin{eqnarray*}
\|f_{n+1}\|_\infty &\le& \|f_n+hG(f_n)\|_\infty\\
&\le& \|f_n\|_\infty + h \left \|\sum_{i=0}^\infty a_i f_n^i \right\|_\infty\\
&\le& \|f_n\|_\infty + h \sum_{i=0}^\infty \|a_i\|_\infty
\|f_n^i\|_\infty\\
&\le& \|f_n\|_\infty + h \sum_{i=0}^\infty \|a_i\|_\infty
\|f_n\|_\infty^i\\
&\le& \|f_n\|_\infty + h g_\infty(\|f_n\|_\infty)
\end{eqnarray*}
Hence the norm of each step of the implicit-explicit iteration will be
controlled by the behavior of the recursion
\begin{equation}
\label{rec_eqn}
f_{n+1}=f_n+hg_\infty(f_n),
\end{equation}
for $f_n,h,a>0$.  Since we are only concerned with short-time
existence and uniqueness, we look specifically at $h=T/N$ and $0\le n
\le N$, for fixed $T>0$ and $N \in \mathbb{N}$.  

\begin{rem}
\label{clever_rem}
The recursion defined by \eqref{rec_eqn} is an Euler solver for 
\begin{equation}
\label{rec_diff_eqn}
\frac{dy}{dt}=g_\infty(y), \text{ with }y(0)=f_0.
\end{equation}
This equation is separable, and $g_\infty$ is analytic near $f_0$, so
there exists a unique solution for the initial value problem
\eqref{rec_diff_eqn} for sufficiently short time.
Also, whenever $y(t)>0$
\begin{equation*}
\frac{d^2 y }{dt^2}=g_\infty'(y(t)) > 0,
\end{equation*}
the function $y(t)$ is concave up.  As a result, the exact solution to
\eqref{rec_diff_eqn} provides an upper bound for the recursion
\eqref{rec_eqn}.  More precisely, we have the following result.
\end{rem}

\begin{lem}
\label{clever_bound}
Suppose $y(0)=f_0 > 0$ in \eqref{rec_diff_eqn}.  Let $T>0$ be
given so that $y$ is continuous on $[0,T]$, and let $N\in\mathbb{N}$.
Then for each $0\le n \le N$, $f_n \le y(T)$, where $f_n$ satisfies
\eqref{rec_eqn} with $h=T/N$.
\begin{proof}
Since the right side of \eqref{rec_diff_eqn} is strictly positive, the
maximum of $y$ is attained at $T$ on any interval $[0,T]$ where $y$ is
continuous.  Furthermore, since $y(0)>0$, it follows from Remark
\ref{clever_rem} that $y$ is concave up on all of $[0,T]$.  Therefore,
$y$ is a convex function on $[0,T]$.  Hence Euler's method,
\eqref{rec_eqn}, will always underestimate the true value of $y$.
Another way of stating this is that
\begin{equation*}
f_n \le y(nh) \le y(T).
\end{equation*}
\end{proof}
\end{lem}


Using Lemma \ref{clever_bound}, the growth of iterates to
\eqref{rec_eqn} may be controlled independently of the step size.
This provides a uniform bound on the sequence of implicit-explicit
approximations.

\begin{lem}
\label{inf_bound}
Suppose $0<h_i=T/i$ for $i \in \mathbb{N}$.  Let $u_i$ be the
implicit-explicit iteration of size $h_i$ beginning at $f_0 \in X$ on
$t \in [0,T]$.  Then there exists a $B>0$ such that for each $i$ and
$t \in [0,T]$, we have $\|u_i(t)\|_\infty \le B$ for sufficiently
small $T>0$.
\begin{proof}
Suppose $f_{in}$ is the $n$-th step of the implicit-explicit iteration
of size $h_i$.  If we let $y(0)=\|f_0\|_\infty$, Lemma
\ref{clever_bound} implies that for any $i$ and any $0 \le n \le i$
\begin{equation*}
\|f_{in}\|_\infty \le y(T)
\end{equation*}
for sufficiently small T.  Hence by \eqref{hubbard_df_fcn} and the
triangle inequality, $\|u_i(t)\|_\infty \le B$ for all $i$ and $t \in
[0,T]$.
\end{proof}
\end{lem}

With the bound on the suprema of the approximations, we can obtain a
bound on the 1-norms.  

\begin{lem}
\label{one_bound}
Suppose $0<h_i=T/i$ for $i \in \mathbb{N}$.  Let $u_i$ be the
implicit-explicit iteration of size $h_i$ beginning at $f_0 \in X$ on
$t \in [0,T]$.  Then there exists an $A>0$ such that for each $i$ and
$t \in [0,T]$, we have $\|u_i(t)\|_1 \le A$ for sufficiently small
$T>0$.
\begin{proof}
First, notice that Lemma \ref{inf_bound} implies that there is a $B>0$
such that for each $i$ and $t \in [0,T]$, we have $\|u_i(t)\|_\infty
\le A$ for sufficiently small $T>0$.  Again suppose $f_{in}$ is the
$n$-th step of the implicit-explicit iteration of size $h_i$.  Then we
compute
\begin{eqnarray*}
\|f_{i,n+1}\|_1 &\le& \|f_{in}\|_1 + h_i \|G(f_{in})\|_1\\
&\le&\|f_{in}\|_1 + h_i \sum_{k=0}^\infty \|a_k f_{in}^k\|_1\\
&\le&\|f_{in}\|_1 + h_i \sum_{k=0}^\infty \int{|a_k f_{in}^k| dx}\\
&\le&\|f_{in}\|_1 + h_i \sum_{k=1}^\infty \|f_{in}\|_\infty^{k-1}
\|a_k\|_\infty \|f_{in}\|_1+h_i\|a_0\|_1\\
&\le&\|f_{in}\|_1\left(1 + h_i \sum_{k=1}^\infty \|a_k\|_\infty
B^{k-1}\right) + h_i\|a_0\|_1\\ 
&\le&\|f_{in}\|_1\left(1 + \frac{h_i}{B} g_\infty(B) -
\frac{h_i}{B}\|a_0\|_\infty \right) + h_i\|a_0\|_1\\ 
&\le&\|f_{in}\|_1\left(1 + h_i C \right) + h_i\|a_0\|_1\\ 
\end{eqnarray*}
This recurence leads to
\begin{eqnarray*}
\|f_{in}\|_1 &\le& \|f_0\|_1 ( 1 + h_i C)^n + h_i \|a_0\|_1
\sum_{m=0}^{n-1} {(1+h_i C)^m}\\
&\le& \|f_0\|_1 ( 1 + h_i C)^n + h_i \|a_0\|_1
\frac{(1+h_i C)^n-1}{h_i C}\\
&\le& \left(\|f_0\|_1 + \frac{1}{C}\|a_0\|_1\right) ( 1 + h_i C)^n -
\frac{1}{C}\|a_0\|_1\\ 
&\le& \left(\|f_0\|_1 + \frac{1}{C}\|a_0\|_1\right) \left( 1 +
\frac{CT}{i} \right)^n - \frac{1}{C} \|a_0\|_1\\
&\le& \left (\|f_0\|_1 + \frac{1}{C}\|a_0\|_1\right) \left( 1 +
\frac{CT}{i} \right)^i - \frac{1}{C} \|a_0\|_1\\ 
&\le&\left( \|f_0\|_1 + \frac{1}{C}\|a_0\|_1 \right) e^{CT} -
\frac{1}{C} \|a_0\|_1 = A. 
\end{eqnarray*}
Once again, by referring to \eqref{hubbard_df_fcn} and using the
triangle inequality, it follows that $\|u_i(t)\|_1 \le B$ for all $i$
and $t \in [0,T]$.
\end{proof}
\end{lem}

\begin{thm}
\label{full_existunique_h_thm}
Suppose $0<h_i=T/i$ for $i\in\mathbb{N}$.  Let $u_i$ be the
implicit-explicit iteration of size $h_i$ beginning at $f_0 \in X$ on
$t \in [0,T]$.  Then, for sufficiently small $T>0$, the sequence
$\{u_i(t)\}_{i=1}^\infty$ converges pointwise to a function in $t$.
The limit function is piecewise differentiable in $t$.
\begin{proof}
This compiles the results of Lemma \ref{existunique_h_thm}, Lemma
\ref{inf_bound}, and Lemma \ref{one_bound}.
\end{proof}
\end{thm}

\begin{rem}
These proofs can be generalized further to handle all equations of
the form

\begin{equation*}
\frac{\partial u(t)}{\partial t} = L(u(t)) + G(u),
\end{equation*}
where $G$ is as in \eqref{def_G}.  If the operator $L$ satisfies
\begin{itemize}
\item $L:L^1(\mathbb{R})\cap L^\infty(\mathbb{R}) \cap
  C^\infty(\mathbb{R}) \rightarrow L^\infty(\mathbb{R}) \cap
  C^\infty(\mathbb{R})$ is a sectorial linear operator \cite{Henry},
\item $\|(I-hL)^{-1}\|_1 \le 1$ and $\|(I-hL)^{-1}\|_\infty \le 1$,
\end{itemize}
then the implicit-explicit iteration
\begin{equation*}
f_{n+1} = (I-hL)^{-1}( f_n + h G(f_n) )
\end{equation*}
converges for whenever $f \in X$.
\end{rem}

\section{Conclusions}

The convergence proof for the implicit-explicit method presented here
has a number of advantages.  First of all, like all implicit-explicit
methods, each approximation to the solution is computed explicitly.
As a result, a fully discretized version (as is standard in the
literature) is easy to program on a computer.  Theorem \ref{full_existunique_h_thm}
therefore assures the convergence of these fully discrete methods.

However, since the implicit-explicit method presented here is
discretized only in time, the convergence proof actually shows the
existence of a semigroup of solutions.  As a result, the convergence
proof forms a bridge between the functional-analytic viewpoint of
differential equations, namely that of semigroups, and the numerical
methods used to approximate solutions.  While the existence and
uniqueness of solutions for \eqref{pde2} has been known via semigroup
methods, the proof provided here gives a more elementary explanation
of how this occurs.  In particular, it approximates the semigroup
action directly.

\bibliography{conv_bib}
\bibliographystyle{plain}

\end{document}